\documentclass[12pt]{article}
\usepackage{amsfonts}
\usepackage{mathrsfs}
\usepackage{amssymb}
\usepackage{amsmath}
\usepackage{latexsym}
\usepackage{graphicx}
\usepackage{color}
 \numberwithin{equation}{section}

\newcommand{\dv}{\mathrm{div}\,}

\newcommand{\p}{\partial}
\newtheorem{Theorem}{Theorem}[section]
\newtheorem{Lemma}{Lemma}[section]
\newtheorem{Proposition}{Proposition}[section]
\newtheorem{Def}{Definition}[section]

\title{\bf A blow-up criterion for compressible viscous
 heat-conductive flows }
\author{
\bf\large Song Jiang \qquad   Yaobin Ou \\
\small LCP, Institute of Applied Physics and Computational Mathematics,\\
\small P. O. Box 8009, Beijing 100088, P.R.China\\
\small E-mail: jiang@iapcm.ac.cn,\quad ou.yaobin@gmail.com}
\begin{document}
\date{}
\maketitle
\begin{abstract}
We study an initial boundary value problem for the Navier-Stokes
equations of compressible viscous heat-conductive fluids in a 2-D periodic
domain or the unit square domain. We establish a blow-up criterion for
the local strong solutions in terms of the gradient of the velocity
only, which coincides with the famous Beale-Kato-Majda criterion for
ideal incompressible flows.




{\bf Running Title:} Blow-up criteria for viscous heat-conductive flows
\end{abstract}
\newpage
\section{Introduction}
This paper is concerned with blow-up criteria for the
two-dimensional Navier-Stokes equations of viscous
heat-conductive gases in a bounded domain $\Omega\subset {\mathbb R}^2$ which
describe the conservation of mass, momentum and total energy, and
can be written in the following form:
\begin{eqnarray}
&& \partial_{t}\rho+\dv(\rho  {u})=0,\label{1.1}\\
&& \rho(\partial_{t}  {u}+ { {u}}\cdot\nabla  { {u}}) -\mu\Delta
 {u}-(\lambda+\mu)\nabla\dv  {u}+\nabla P=0,\ \label{1.2}\\
&& c_{_V}\rho ( \partial_{t}\theta+ {u}\cdot\nabla\theta)
-\kappa\Delta\theta +P\dv  {u}=\frac{\mu}{2}|\nabla
 { {u}}+\nabla  {u}^{t}|^{2}+\lambda(\dv
 {u})^{2}. \label{1.3}
\end{eqnarray}
Here we denote by $\rho, \theta$ and $ {u}=(u_1,u_2)^t$ the
density, temperature, and velocity, respectively. The physical
constants $\mu, \lambda$ are the viscosity coefficients satisfying
$\mu>0$, $\lambda + \mu \geq 0$, $c_{V}>0$ and $\kappa >0$ are the
specific heat at constant volume and thermal conductivity
coefficient, respectively. $P$ is the pressure which is a known
function of $\rho$ and $\theta$, and in the case of an ideal gas $P$
has the following form
\begin{equation}\label{es}
P=R\rho\theta,
\end{equation}
where $R>0$ is a generic gas constant.

Let $\Omega$ be a periodic domain in $\mathbb{R}^2$, or the unit square
$[0,1]^2$ in $\mathbb{R}^2$. We will consider an initial boundary
value problem for \eqref{1.1}--\eqref{1.3} in
$Q:=(0,\infty)\times\Omega$ with initial condition
\begin{equation}\label{1.4}
 (\rho,u,\theta)|_{t=0}
=(\rho_{0},u_{0},\theta_{0})\ \ {\rm in}\ \ \Omega,
\end{equation}
and boundary conditions:
\begin{equation}\label{1.5}
\begin{array}{l}
 (i)\qquad\rho,\,u,\,\theta\; \textrm{are\; periodic\; in\; each}\; x_i\;
\textrm{\textrm{for}}\; 1\le i\le 2,\; \textrm{or}\\[2mm]
(ii)\qquad u_i|_{x_i=0,1}=\partial_{i}u_j|_{x_i=0,1}=0,\;
\forall\; 1\le i,j\le 2,\,j\neq i,\\
\qquad\quad \partial_i\theta|_{x_i=0,1}=0,\; \forall\;1\le i\le 2.
\end{array}
\end{equation}

In the last decades significant progress has been made in the study
of global in time existence for the system (\ref{1.1})--(\ref{1.5}).
With the assumption that the initial data are sufficiently small,
Matsumura and Nishida \cite{9,10} first proved the global existence
of smooth solutions to initial boundary value problems and the
Cauchy problem for \eqref{1.1}--\eqref{1.3}, and the existence of
global weak solutions was shown by Hoff \cite{Hoff}. For large data,
however, the global existence to
(\ref{1.1})--(\ref{1.5}) is still an open problem, except certain special
cases, such as the spherically symmetric case in domains without the
origin, see \cite{Ji96} for example. Recently, Feireisl
\cite{Fe04a,Fe04b} obtained the global existence of the so-called
``variational solutions'' to (\ref{1.1})--(\ref{1.3}) in the case of
real gases in the sense that the energy equation is replaced by an
energy inequality. However, this result excludes the case of ideal
gases unfortunately. We mention that in the isentropic case, the
existence of global weak solutions of the multidimensional
compressible Navier-Stokes equations was first shown by Lions
\cite{Lions98}, and his result was then improved and generalized in
\cite{FNP01} (also see \cite{JZ00,JZ03}, and among others). Moreover, this kind of
weak solution with finite energy was shown to exist in $[0,\infty)$
as long as the density remains bounded in $L^\infty(\mathbb{T}^2)$
(cf. Desjardins \cite{DES}).

Xin \cite{22}, Rozanova \cite{Roz08} showed the non-existence of
global smooth solutions when the initial density is compactly
supported, or decreases to zero rapidly. Since the system
(\ref{1.1})--(\ref{1.3}) is a model of non-dilute fluids, these
non-existence results are natural to expect when vacuum regions are
present initially. Thus, it is very interesting to investigate
whether a strong or smooth solution will still blow up in finite
time, when there is no vacuum initially. Recently, Fan and Jiang
\cite{3} proved the following blow-up criteria for the local strong
solutions to (\ref{1.1})--(\ref{1.5}) in the case of two dimensions:
$$
\lim_{T\to T^*}\Big( \sup_{0\leq t\leq
T}\{\|\rho\|_{L^\infty},\|\rho^{-1}\|_{L^\infty},
\|\theta\|_{L^\infty}\}(t)
+\int^T_0(\|\rho\|_{W^{1,q_0}}+\|\nabla\rho\|^4_{L^2}
+\|u\|^{\frac{2r}{r-2}}_{L^{r,\infty}})dt \Big)=\infty,
$$
or,
$$
\lim_{T\to T^*} \Big ( \sup_{0\leq t\leq
T}\{\|\rho\|_{L^\infty},\|\rho^{-1}\|_{L^\infty},
\|\theta\|_{L^\infty}\}(t)
+\int^T_0(\|\rho\|_{W^{1,q_0}}+\|\nabla\rho\|^4_{L^2})dt
\Big)=\infty , $$
provided $2\mu>\lambda$, where $T^*<\infty$ is the maximal time of
existence of a strong solution $(\rho,u)$, $q_0>3$ is a certain
number, $3<r\leq\infty$ with $2/s+3/r=1$, and $L^{r,\infty}\equiv
L^{r,\infty}(\Omega)$ is the Lorentz space.

In the isentropic case, the result in \cite{3} reduces to
\begin{equation}\label{fj}
\lim\limits_{T\rightarrow T^*}\Big(\sup\limits_{0\le t\le
T}\|\rho\|_{L^\infty}+
\int_0^T\left(\|\rho\|_{W^{1,q_0}}+\|\nabla\rho\|_{L^2}^4\right)\Big)=\infty,\quad
\end{equation}
provided $7\mu> 9\lambda.$ Recently, Huang and Xin \cite{HX09}
established the following blow-up criterion in a  3-D smooth bounded
domain, similar to the Beale-Kato-Majda criterion for ideal
incompressible flows \cite{BKM},
 for the isentropic compressible Navier-Stokes equations:
\begin{equation}\label{bu}
\lim\limits_{T\rightarrow T^*}\int_0^T\|\nabla
u\|_{L^\infty}dt=\infty,
\end{equation}
provided
\begin{equation}
7\mu > \lambda. \label{coef}
\end{equation}
Indeed, if the domain is a periodic or unit square domain in
$\mathbb{R}^2$, the blow-up criterion is refined by Fan, Jiang and
Ni \cite{FJN},  to be
\begin{equation}
\lim\limits_{T\rightarrow T^*}\sup\limits_{0\le t\le T}
(\|\rho\|_{L^\infty}+\| {\rho}^{-1}\|_{L^\infty})=\infty.
\end{equation}
This result was recently improved by Sun, Wang and Zhang \cite{SZ10},
in both two- and three-dimensional cases, to be
$$
\lim\limits_{T\rightarrow T^*}\sup\limits_{0\le t\le T}\|\rho\|_{L^\infty}=\infty ,
$$
while a sharper criterion in the following form was given by Haspot in \cite{Ha10}:
$$
\lim\limits_{T\rightarrow T^*}\sup\limits_{0\le t\le T}\|\rho\|_{L^{(N+1+\epsilon)\gamma}}=\infty ,
$$
where $N$ $(=2,3)$ and $\gamma$ are the spatial dimension and the specific heat ratio respectively,
and $\epsilon$ is an arbitrary small number.

For the non-isentropic compressible Navier-Stokes equations, Fan, Jiang and Ou
\cite{FJO} established a blow-up criterion with additional upper
bound of $\theta$:
$$\lim_{T\rightarrow T^*}\Big( \|\theta\|_{L^\infty (0,T;L^\infty )}
+\|\nabla u\|_{L^1(0,T;L^\infty )}\Big) =\infty ,$$ provided that
the condition (\ref{coef}) is satisfied. This result reduces to the
one in \cite{HX09} in the isentropic regime.

The aim of this paper is to show that, the requirement of upper boundedness
of $\theta$ in \cite{FJO} can be removed and the condition
\eqref{coef} can be refined to be the usual physical condition
$\lambda+  \mu \ge 0$ for non-vacuum fluids in a 2-D domain.  This
result coincides the famous Beale-Kato-Majda criterion for ideal
incompressible flows, and the criterion in \cite {HX09} in the
non-vacuum case. In contrast to \cite{FJO}, it is interesting to see here
that the temperature $\theta$ allows to vanish in $\Omega$, and more
important, the temperature will not lead to the blow-up of strong solutions
to the full Navier-Stokes equations. These are exactly the new points of this paper,
in comparison with \cite{FJO}.

 Moreover, it is interesting to
see that the a priori assumption \eqref{2.1} is more concise than
the one in \cite{3,FJO}. 


For the sake of generality, we will study the blow-up criterion for
local strong solutions.

Before giving our main result, we state the following local
existence of the strong solutions, the proof of which can be found
in \cite{21}.

\begin{Proposition}\label{pr1.1}{\rm (Local Existence)}
Let $\Omega$ be a bounded domain in $\mathbb{R}^2$ as previously
stated. Suppose that the initial data $\rho_{0},u_{0},\theta_{0}$
satisfy
\begin{equation}
\begin{split}
& \inf\limits_{x\in\Omega}{ \rho_0 (x)> 0,}\quad \rho_{0}\in W^{1,q}
(\Omega)\quad\mbox{for any $q>2$},   \\
& u_{0}\in H_{0}^{1}(\Omega)\cap H^{2}(\Omega),\;\; \theta_0\ge
0,\;\; \theta_{0}\in H^{2}(\Omega),
\end{split}
\label{1.6}
\end{equation}
and the compatibility conditions
\begin{equation}\label{compat}
\begin{split}
&\mu\Delta u_0 +(\mu+\lambda)\nabla{\rm div}u_0-
R\nabla(\rho_0\theta_0)=\rho_0^{1/2}g_1,
\\
&\kappa\Delta \theta_0 +\frac{\mu}{2}|\nabla u_0+\nabla
u^{t}_0|^{2}+\lambda(\dv u_0)^{2}-R\rho_0\theta_0 {\rm
div}u_0=\rho_0^{1/2}g_2,
\end{split}
\end{equation}
$\mbox{for some}\;g_1,g_2\in L^2(\Omega ).$
 Then there exist a positive constant $T_0$ and a unique strong solution
$(\rho,u,\theta)$ to (\ref{1.1})--(\ref{1.5}), such that
\begin{equation}
\begin{array}{rl} &
{  \rho > 0}, \;\rho\in C([0,T_0];W^{1,q}),\;
\rho_{t}\in C([0,T_0];L^{q}),\\[1mm]
& u\in C([0,T_0];H_{0}^{1}\cap H^{2})\cap
L^{2}(0,T_0;W^{2,q}),\;\;\\
& u_{t}\in L^{\infty}(0,T_0;L^{2}),\;u_t\in L^{2}(0,T_0;H_{0}^{1}),\\[1mm]
& {\theta\ge 0},\; \theta\in C([0,T_0];H^{2})\cap
L^{2}(0,T_0;W^{2,q}),\\
& \theta_{t}\in L^{\infty}(0,T_0;L^{2}),\;\theta_t\in
L^{2}(0,T_0;H^{1}).
\end{array}\label{1.7}
\end{equation}
\end{Proposition}

By the regularity $(u_t,\theta_t)\in L^{\infty}(0,T_0;L^2)$,
the quantities $\|u_t(T_0)\|_{L^2(\Omega)}$
and $\|\theta_t(T_0)\|_{L^2(\Omega)}$, redefined if necessary, are
finite, which leads to the validity of the compatibility conditions
at $t=T_0$. One may refer to \textit{Remark 2} in \cite{21} for the
necessity of the compatibility conditions in \eqref{compat}.

Therefore, with the regularities in \eqref{1.7} and the new
compatibility conditions at $t=T_0$, we are able to extend the
solution to the time beyond $T_0$. Now, we are interested in the
question what happens to the solution if we extend the solution
repeatedly. One possible case is that the solution exists in
$[0,\infty )$, while another case is that the solution will blow-up
in finite time in the sense of \eqref{1.7}, that is, some of the
regularities in \eqref{1.7} no longer hold.

{\Def $T^*\in (0,\infty)$ is called the maximal life-time of
existence of a strong solution to \eqref{1.1}-\eqref{1.5} in the
regularity class \eqref{1.7} if for any $0<T<T^*$, $(\rho,u,\theta)$
solves \eqref{1.1}-\eqref{1.5} in $[0,T]\times\Omega$ and satisfies
\eqref{1.7} with $T_0=T$, and moreover, \eqref{1.7} does not hold for
$T_0=T^*$.}

\vskip .2cm
  Now, we are ready to state the main theorem of this paper.
\begin{Theorem}\label{th1.1} {\rm (Blow-up Criterion)} Suppose that
the assumptions in Proposition 1.1 are satisfied. Let
$(\rho,u,\theta)$ be the strong solution obtained in Proposition
1.1. Then either this solution can be extended to $[0,\infty)$, or
there exists a positive constant $T^*<\infty$, the maximal time of
existence, such that the solution only exists in $[0,T]$ for every
$T<T^*$, and
$$\lim_{T\rightarrow T^*}\int_0^T \|\nabla u(t)\|_{L^\infty}dt=\infty.$$
\hfill $\Box$
\end{Theorem}

We will prove Theorem \ref{th1.1} by contradiction in the next
section. In fact, the proof of the theorem is based on a priori
estimates under the assumption that $\|\nabla
u\|_{L^1(0,T;L^\infty)}$ is bounded independent of any $T\in
[0,T^*)$. The a priori estimates are then sufficient for us to apply
the local existence theorem repeatedly to extend the local solution
beyond the maximal time of existence $T^*$, consequently,
contradicting the maximality of $T^*$.

The proof of this paper is based on the estimates for the effective
viscous flux $(2\mu+\lambda)\dv u-P$. We can obtain   good estimates
on the effective viscous flux at the first step of derivative
estimates. This is the main ingredient of the estimates. It plays an
important role in deriving other derivative estimates. This
technique is applied in many previous situations for studying the
Navier-Stokes equations (cf. \cite{VK,FJN}), we adapt it here to
establish the blow-up criteria for the full Navier-Stokes
equations.

The rest of this paper is organized as follows.
First, we will establish the estimates for all the zero-th order quantities of the solutions.
Then, we will derive the crucial $\|u\|_{C([0,T],H^1(\Omega))}$ bound
by utilizing the effective viscous flux, and then the estimates for other derivatives.
Finally, we conclude the blow-up criteria by contradiction and continuity arguments.

Throughout this paper, we will use the following abbreviations:
$$L^p\equiv L^p(\Omega ),\quad H^m\equiv H^m(\Omega ), \quad H^m_0\equiv H^m_0(\Omega ), $$
$$L^p_t(X)\equiv L^p(0,t;X(\Omega)),\quad  C_t(X)\equiv C([0,t],X(\Omega)).$$

\section{Proof of Theorem \ref{th1.1}}
Let $0<T<T^*$ be arbitrary but fixed. Throughout this section, We
 denote by $\delta,\epsilon$ various small positive constants, and
 moreover
we denote by $C$ (or $C_X$ to emphasize the dependence of $C$ on $X$
) a general positive constant which may depend continuously on
$T^*$.

Let $(\rho,u,\theta )$ be a strong solution to the problem
(\ref{1.1})--(\ref{1.5}) in the function space given in (\ref{1.7})
on the time interval $[0,T]$. Suppose that $T^*<\infty$. We will
prove Theorem \ref{th1.1} by a contradiction argument. To this end,
we suppose that for any $T<T^*$,
\begin{equation}
\|\nabla u\|_{L^1(0,T;L^\infty )}\leq C<\infty, \label{2.1}
\end{equation}
we will deduce a contradiction to the maximality of $T^*$.

\subsection{Zero-th order estimates}
First, we show that the density $\rho$ is bounded from below and
above due to the assumption  in (\ref{2.1}). It is easy to see that
the continuity equation (\ref{1.1}) on the characteristic curve
$\dot{\chi}(t)=u(\chi(t))$ can be written as
$$\frac{d}{dt}\rho(\chi(t),t)= -\rho (\chi(t),t) {\rm div}\,u (\chi(t),t).$$
Thus, by Gronwall's inequality and (\ref{2.1}), one obtains that for
any $x\in\bar\Omega$ and $t\in [0,T]$,
\begin{equation}
\begin{split}
 C^{-1}\leq\underline{\rho}\exp(-\|{\rm
 div}\,u\|_{L^1_T( L^\infty)})
 \leq \rho (x,t)\le \bar{\rho}\exp(\|{\rm
 div}\,u\|_{L^1_T( L^\infty)})\leq C, \label{2.1a}
 \end{split}
 \end{equation}
where $0< \underline{\rho}\le \rho_0\le \bar{\rho}$.

Multiplying \eqref{1.2} by $u$ and summarizing the result by \eqref{1.3}, we
integrate to get
\begin{equation}\label{s5}
\begin{split}
\frac{d}{dt}\int_\Omega\Big( \frac{1}{2}\rho |u|^2+ c_V\rho\theta\Big)dx=0.
 \end{split}
 \end{equation}

 Next, we show
that $\theta$ is non-negative in $[0,T]\times\Omega$ (see also
\cite{Fe04a,FJO}). Let $H(\theta)(x,t):=c_V\min\{-\theta(x,t),0\}$.
Clearly, $H'(\theta)\le 0$ and $H''(\theta)=0$. We multiply
\eqref{1.3} by $H'(\theta)$ and integrate over $\Omega$ to obtain
\begin{equation*}
\begin{split}
\int_\Omega\Big(\rho& (H(\theta)_t+ u\cdot\nabla H(\theta))+R\rho
H(\theta)\dv u\Big)dx\\
&=\int_\Omega  H'(\theta)\Big( \kappa\Delta \theta +
\frac{\mu}{2}|\nabla u+\nabla u^{t}|^{2}+\lambda(\dv u)^{2}\Big)
dx\\
&\le \kappa\int_{\partial\Omega}H'(\theta)\frac{\partial\theta}{\partial n}dS
-\kappa\int_\Omega H''(\theta)|\nabla \theta|^2dx \\
& \le 0
 \end{split}
 \end{equation*}
 By the continuity equation \eqref{1.1}, we integrate by parts to
 get
\begin{equation*}
\begin{split}
\frac{d}{dt}\int_\Omega \rho  H(\theta) dx  &\le C\int_\Omega |\dv
u||\rho
H(\theta)|dx\\
&\le -\|\dv u\|_{L^\infty}\int_\Omega  \rho H(\theta) dx.
 \end{split}
 \end{equation*}
 Utilizing \eqref{2.1} and Gronwall's inequality, we have
 $$\int_\Omega \rho H(\theta) dx\equiv 0,\quad \forall\; t\in [0,T],$$
 since $\theta_0\ge 0$. Thus $\theta\ge 0$ by the definition of $H(\theta)$.

From \eqref{2.1a}, \eqref{s5} and the non-negativeness of $\theta$,
we have
\begin{equation}
\sup\limits_{0\le t\le T}(\|u(t)\|_{L^2}+\|\theta(t)\|_{L^1})\le C,\quad \forall\;t\in [0,T].
\end{equation}

\subsection{Estimates for derivatives}

By multiplying \eqref{1.2} by $u$ and integrating by parts, it follows immediately 
\begin{equation}
\begin{split}
\|\nabla u   \|_{L^2_t( L^2)}^2  &\le
\frac{1}{2}\int_\Omega\rho|u|^2 dx (t) +
\frac{1}{2}\int_\Omega \rho_0 |u_0|^2 dx + \int_0^t\int_\Omega|P||\dv u|dxds\\
&\le C+R\|\rho\|_{L^\infty_{t }(L^\infty)}\|\theta\|_{L^\infty_t(
L^1)}\|\dv u\|_{L^1_t( L^\infty)} \le C,\;\forall \;t\in[0,T].
\end{split}
\end{equation}

Now, we are ready to control $\|u\|_{L^\infty (0,t;H^1)}$
and $\|u\|_{L^2(0,t;H^2)}$ by estimating the effective viscous flux
$(2\mu + \lambda)\dv u-P$, which is similar to the strategies in \cite{VK} and \cite{FJN}.
These are the key estimates in our proof. To simplify the statement, we denote by
$$V:=-\textrm{curl} u=\partial_2 u_1-\partial_1 u_2$$
the vorticity, and by
$$F:=(2\mu + \lambda)\dv u - P $$
the effective viscous flux.

\begin{Lemma}{\rm (Key estimates)}\label{eff} Under the assumption in \eqref{2.1},
we have for any $T<T^*$ that
\begin{eqnarray} &&
\sup\limits_{0\le t\le T}\|(\theta,\nabla u)(t)\|_{L^2}^2
+\int_{0}^{T} \|(\nabla\theta, \nabla F)\|_{L^2}^2dt\leq C.
\label{2.3}
\end{eqnarray}
\end{Lemma}

\noindent\textbf{Proof.} We first derive the following system from
\eqref{1.1}-\eqref{1.2} and \eqref{es}-\eqref{1.5}:
\begin{eqnarray}
&& V_t+ u\cdot\nabla V+ V\dv u- \partial_2\Big(\frac{1}{\rho}(\partial_2V+\p_1 F)\Big)
    + \p_1 \Big(\frac{1}{\rho}(\p_2 F-\p_1 V )\Big)=0,\label{e1}\\
    && F_t+ u\cdot\nabla F-(2\mu+\lambda) \Big(\partial_1\Big(\frac{1}{\rho}(\partial_2V+\p_1 F)\Big)
    + \p_2 \Big(\frac{1}{\rho}(\p_2 F-\p_1 V )\Big)\Big)\nonumber\\
    &&\qquad\qquad\quad= O(1)(\p_i u_j \p_k u_l)-\gamma P\dv u-(\gamma-1)\Delta\theta,
    \label{e2}\\  [2mm]
&& V|_{t=0}={\rm curl}u_0,\; F|_{t=0}= (2\mu+\lambda)\dv u_0-R\rho_0\theta_0,\label{e3}\\
&& V|_{x_i=0,1}=\p_i F|_{x_i=0,1}=0,\; i=1,2,\label{e4}
\end{eqnarray}
where $\gamma=1+R/c_V$. We multiply \eqref{e1}, \eqref{e2} by $V$,
$F/(2\mu +\lambda)$ respectively, and integrate by parts to get
\begin{equation}\label{1.10}
\begin{split}
\frac{1}{2}\frac{d}{dt}& \int_\Omega(V^2 +
\frac{F^2}{2\mu+\lambda})dx
+\int_\Omega\frac{(\p_2 V+ \p_1 F)^2 + (\p_2 F-\p_1 V)^2}{\rho} dx\\
& =\int_\Omega \frac{1}{2}\dv u(\frac{F^2}{2\mu
+\lambda}-V^2)-\frac{\gamma}{2\mu+\lambda}
\int_\Omega FP\dv u dx\\
&\quad +O(1)\int_\Omega \p_iu_j\p_ku_l F
dx-\frac{\gamma-1}{2\mu+\lambda}\int_\Omega F\Delta \theta dx \\
& =:\sum_{i=1}^4 I_i.
\end{split}
\end{equation}
Denote the second integral on the left-hand side by $I_5$.
Then by integrating by parts and applying \eqref{2.1a}, we have
\begin{equation}
\begin{split}
I_5 &\ge C_0\int_\Omega \big[ (\p_2 V+\p_1 F)^2 +(\p_2 F- \p_1 V)^2\big] dx\\
&= C_0 (\|\nabla V\|_{L^2}^2 +\|\nabla F\|_{L^2}^2).
\end{split}
\end{equation}
Noting that $\rho$ is bounded from above, we obtain
\begin{equation}
\begin{split}
|I_1+I_2| &\le C\|\dv u\|_{L^\infty}(\|V\|_{L^2}^2 +\|F\|_{L^2}^2
+\|\theta\|_{L^2}^2).
\end{split}
\end{equation}

To bound $I_3$, we need the following lemma.
\begin{Lemma}
 For any $u\in H^1(\Omega)$ satisfying the boundary conditions
in \eqref{1.5}, we have
\begin{equation*}
\|\nabla u\|_{L^p}\le C(\|\dv u\|_{L^p} + \|{\rm curl} u\|_{L^p} +
\|u\|_{L^2}),\quad\forall\; 2\le p <\infty.
\end{equation*}
\end{Lemma}
The previous version of this lemma (cf. \cite{BB}) holds in case of a
smooth domain $\Omega$
 and $u\cdot n|_{\p\Omega}=0$. However, the conclusion can be easily adapted to our case.
 We can slightly modify the original proof by an
 extension argument, since the angles at corner points of our domain here are right-angles.

 Noting also that
\begin{equation*}
\dv u =(F+P)/(2\mu +\lambda),
\end{equation*}
we have
\begin{equation*}
\begin{split}
|I_3| &\le C\|\nabla u\|_{L^\infty}\|\nabla u\|_{L^2}\|F\|_{L^2}\\
 &\le C\|\nabla u\|_{L^\infty}(\|V\|_{L^2}+\|F\|_{L^2}+\|\theta\|_{L^2}+ \|u\|_{L^2})\|F\|_{L^2}\\
 &\le C\|\nabla u\|_{L^\infty}(\|V\|_{L^2}^2+\|F\|_{L^2}^2+\|\theta\|_{L^2}^2+ \|u\|_{L^2}^2).
\end{split}
\end{equation*}
Finally, we get by integration by parts that
\begin{equation*}
\begin{split}
|I_4| &\le C\int_\Omega |\nabla\theta||\nabla F|dx
\le \frac{C_0}{2}\|\nabla F\|_{L^2}^2 + C_1\|\nabla \theta\|_{L^2}^2.
\end{split}
\end{equation*}
On the other hand, we can derive from \eqref{1.3} that
\begin{equation}\label{1.13}
\begin{split}
\frac{c_V}{2}\frac{d}{dt} \int_\Omega\rho \theta^2dx
+\int_\Omega|\nabla\theta|^2 dx & \le C(\|\dv u\|_{L^\infty}\|\theta\|_{L^2}^2
 +\|\theta\|_{L^2}\|\nabla u\|_{L^\infty}\|\nabla u\|_{L^2})\\[1mm]
 &\le C \|\nabla u\|_{L^\infty}(\|V\|_{L^2}^2+\|F\|_{L^2}^2+\|\theta\|_{L^2}^2).
\end{split}
\end{equation}

Collecting the above estimates,
 we can draw the conclusion  by applying
 Gronwall's inequality to \eqref{1.10} and \eqref{1.13}.

\hfill $\Box$ \vskip .2cm

\begin{Lemma}\label{lem2.5}
With the assumption in \eqref{2.1}, we have for any $0\le t\le T$ that
$$ \|\nabla\rho (t)\|_{L^2}+  \|u\|_{L^2_t(H^2)} + \|u_t\|_{L^2_t(L^2)}
 \le C (1+\|\theta\|_{L^2_t(H^2)}^{4\epsilon}) .$$
\end{Lemma}
\textbf{Proof.} Since $u$ is a solution of the strictly elliptic system
$$-\mu\Delta u=f$$
where $f:=-\rho u_{t}-\rho u\cdot\nabla u-\nabla F$, it follows from the classical
regularity theory and the interpolation inequality that
\begin{eqnarray}\label{uh20}
\|u\|_{H^{2}} &\leq & C\left(\|
u_{t}\|_{L^{2}}+\|u\|_{H^1}^\frac{3}{2}\|\nabla
u\|_{H^2}^\frac{1}{2} +\|\theta\|_{L^\infty}\|\nabla
\rho\|_{L^{2}}+\|\nabla\theta\|_{L^2} \right),
\end{eqnarray}
whence,
\begin{equation}\label{ul1h2}
\|u\|_{L^1_t( H^{2})} \leq C(1+\| u_{t}\|_{L^{2}_t( L^2)}
+\int_0^t\|\theta\|_{L^\infty}\|\nabla \rho\|_{L^{2}}ds),\quad
\forall\; t\in [0,T].
\end{equation}

On the other hand, we would like to estimate $\| u_{t}\|_{L^{2}(0,t;L^2)}$
in terms of $\|u\|_{L^2(0,t;H^{2})}$ to close the estimates.
Taking the inner product of \eqref{1.2} with $u_t$ in $L^2((0,t)\times\Omega)$, we obtain
\begin{equation}
\begin{split}
& \frac{d}{dt}\int_\Omega\Big(\frac{\mu}{2}|D (u)|^{2}
+\frac{\lambda}{2}(\dv u)^{2}-P {\rm div}\,u\Big)dx
+\int_\Omega\rho u_{t}^{2}dx \\[1mm]
& \qquad \leq C\int_\Omega |u||\nabla u||u_t| dx + \int_\Omega P_t
\dv u dx:= J_1 + J_2, \label{aa}
\end{split}
\end{equation}
where $D(u):=(\nabla u + \nabla u_t)/2$. We calculate $J_1$ and $J_2$
as follows. By the interpolation inequality again, we get
\begin{equation*}
\begin{split}
\int_0^t J_1 ds\le & C\int_0^t\|u\|_{H^1}^\frac{3}{2}\|u\|_{H^2}^\frac{1}{2}\|u_t\|_{L^2}ds\\
\le & \frac{1}{2}\int_0^t \int_\Omega\rho u_{t}^{2}dxds+ C
\|u\|_{L^2_t( H^{2})}.
\end{split}
\end{equation*}

From \eqref{1.1} and \eqref{1.3}, we have
\begin{equation}
P_t+u\cdot\nabla P + \gamma P\dv u-\kappa\Delta\theta= 2\mu |D(u)|^2
+ \lambda (\dv u)^2,
\end{equation}
thus, by virtue of integration by parts,
\begin{equation*}
\begin{split}
\Big|\int_0^t J_2 ds\Big|\le &  C\int_0^t\int_\Omega (P(|\nabla u|^2
+ |u||\nabla^2 u|) + |\nabla \theta||\nabla \dv u|
 + |\nabla u|^3)dxds\\
 \le & C\|\theta\|_{L^\infty_t(L^2)}\|\nabla u\|_{L^1_t(L^\infty)}\|\nabla u\|_{L^\infty_t(L^2)}\\
    &+ C\|\theta\|_{L^\infty_t(L^2)}^\frac{1}{2}\|\theta\|_{L^2_t(H^1)}^\frac{1}{2} \|u\|_{L^2_t(H^2)}\|u\|_{L^\infty_t(H^1)}\\
    &+ C\|\nabla \theta\|_{L^2_t(L^2)}\|u\|_{L^2_t(H^2)}+ C\|\nabla u\|_{L^1_t(L^\infty)}\|\nabla u\|_{L^\infty_t(L^2)}^2\\
    \le & C + C\|u\|_{L^2_t(H^2)}.
\end{split}
\end{equation*}
Note that $\int_\Omega P {\rm div}udx (t)\le
C\|\theta\|_{C_t(L^2)}\|\nabla u\|_{C_t(L^2)}\le C,$ by the previous estimates.
Thus, we conclude
\begin{equation}\label{utl2}
\|u_t\|_{L^2_t(L^2)}\le C (1+  \|u\|_{L^2_t(H^2)})^\frac{1}{2},
\quad \forall \; t\in [0,T].
\end{equation}

Next, we apply $\nabla$ to the equation \eqref{1.1}, then multiply
the resulting equation by $\nabla \rho$ and integrate over $\Omega$
to get
\begin{equation*}
\begin{split}
\frac{d}{dt}\int_\Omega |\nabla\rho|^2 dx
 \leq C\|\nabla
u\|_{L^\infty} \|\nabla\rho\|^2_{L^2} +
C\|u\|_{H^2}\|\nabla\rho\|_{L^2}, \label{aa1}
\end{split}
\end{equation*}
which gives,   by Gronwall's inequality,
\begin{equation}
 \|\rho (t)\|_{L^2}\le C (1 +\|u\|_{L^1_t(H^2)}).
\end{equation}
Substituting \eqref{ul1h2} and \eqref{utl2} into the above
inequality and applying the integro-type Gronwall inequality, we
have
\begin{equation}
 \|\rho (t)\|_{L^2}\le C (1 +\|u\|_{L^2_t(H^2)}^\frac{1}{2}).
\end{equation}

From \eqref{uh20} again, we derive
\begin{equation}
\begin{split}
\|u\|_{L^2_t( H^{2})} &\leq C(1+\| u_{t}\|_{L^2_t( L^2)} +
\|\theta\|_{L^2_t(L^\infty)}\|\nabla
\rho\|_{L^\infty_t(L^{2})})\\
& \le C + C \|\theta\|_{L^2_t(H^1)}^{2(1-\epsilon)}
\|\theta\|_{L^2_t(H^2)}^{2 \epsilon } (1
+\|u\|_{L^2_t(H^2)}^\frac{1}{2}),\quad \forall\; t\in [0,T],
\end{split}
\end{equation}
which gives this lemma immediately.\hfill $\Box$

 \vspace{2mm} Next, we derive further estimates for the
derivatives of $\theta$ to close the above estimates.
\begin{Lemma}\label{lemtheta}
Assuming \eqref{2.1},
we have for any $T<T^*$ that
$$
\sup\limits_{0\le t\le T}\|\nabla \theta(t)\|_{L^2}^2 + \int_0^T
(\|\theta\|_{H^2}^2 + \|\theta_t\|_{L^2}^2)dt\le C.$$
\end{Lemma}
\textbf{Proof.}  Multiplying \eqref{1.3} by $\theta_t$ and integrating over $\Omega$, we get
\begin{equation}
\begin{split}
 \frac{\kappa}{2}\frac{d}{dt}& \int_\Omega  |\nabla\theta|^2dx+c_{_V}\int_\Omega\rho\theta_t^2dx\\
 \le & C\int_\Omega (|\rho||u||\nabla \theta| + |\rho||\theta||\dv u|
 + |\nabla u|^2)|\theta_t|dx \\[1mm]
 \le &  C(\|u\|_{H^1}\|\theta\|_{H^1}^{1-\epsilon}\|\theta\|_{H^2}^\epsilon
  + \|\theta\|_{L^\infty}\|\dv u\|_{L^2} +
    \|\nabla u\|_{L^4}^2)\|\theta_t\|_{L^2} \\[2mm]
 \le & \delta \|\theta_t\|_{L^2}^2 +
  C_\delta\|u\|_{H^1}^2 (\|\theta\|_{H^1}^{2-2\epsilon}\|\theta\|_{H^2}^{2\epsilon}
    +  \|u\|_{H^2}^2).
\end{split}
\end{equation}
 Since $\rho$ is bounded from below, we choose $\delta$ small enough
and then apply Gronwall's inequality to conclude
\begin{equation*}
\|\nabla \theta(t)\|_{L^2}^2 + \|\theta_t\|_{L^2_t(L^2)}^2\le C (1+
\|\theta\|_{L^2_t(H^2)}^{4\epsilon}),\quad\forall\;0\le t\le T.
\end{equation*}

Again, we apply the elliptic regularity theory to derive
\begin{equation*}
\begin{split}
 \|\theta\|_{L^2_t(H^2)}& \le C(\|\theta_t\|_{L^2_t(L^2)} + \|u\cdot\nabla \theta\|_{L^2_t(L^2)}
    + \|\theta \dv u\|_{L^2_t(L^2)} + \|\nabla u\|_{L^2_t(L^4)}^2 + \|\theta\|_{L^2_t(H^1)})\\
    &\le  C(\|\theta_t\|_{L^2_t(L^2)} +  \|u\|_{L^2_t(H^2)} \|\theta \|_{L^\infty_t(H^1)} +
     \|u\|_{L^\infty_t(H^1)}  \|u\|_{L^2_t(H^2)}   + 1)\\
     &\le C  (\|\theta_t\|_{L^2_t(L^2)} + \|\theta \|_{L^\infty_t(H^1)} +1).
\end{split}
\end{equation*}
Thus we can show this lemma easily from the above two inequalities
by choosing $\epsilon <1/4$.\hfill $\Box$ \vskip .2cm

Next, we will exploit the a priori estimates obtained so far to derive
bounds on temporal derivatives and high-order derivatives.
\begin{Lemma}\label{lemt} Let  \eqref{2.1} hold. For any $t\in [0,T]$, we have
\begin{equation*}
\sup\limits_{0\le t\le T} \Big(\|(u_t(t), \theta_t(t))\|_{L^2}^2 +\|(u(t),\theta(t))\|_{H^2}^2\Big)
+ \int_0^T \|(u_t,\theta_t) \|_{H^1}^2ds\le C.
\end{equation*}

\end{Lemma}

\noindent{\bf Proof.} Now, taking $\partial_t$ to the equation (\ref{1.2}),
multiplying then the resulting equation by $u_t$ in $L^2(\Omega )$, integrating by parts, and
employing \eqref{1.1} and (\ref{2.3}), we find
\begin{eqnarray}
&& \frac{1}{2}\frac{d}{dt}\int_\Omega\rho
u_t^2dx+\int_\Omega\Big(2\mu|D(u_t)|^2+ \lambda (\dv u_t)^2\Big)dx\nonumber\\
&& \quad =\int_\Omega P_t\dv u_t dx-\int_\Omega\rho u\cdot\nabla\big[
(u_t+u\cdot\nabla u)u_t\big]dx -\int_\Omega\rho u_t\cdot\nabla u\cdot u_t dx
\nonumber\\
&& \quad := K_1+K_2+K_3. \label{2.17}
\end{eqnarray}

Observing that $P_t=R\rho_t\theta+R\rho\theta_t$ and $\rho_t=-\dv
(\rho u)$, and applying the interpolation inequality in two dimensions, we deduce
\begin{equation*}
\begin{split}
|K_1| &\le \delta \|\nabla u_t\|_{L^2}^2+ C_\delta\Big(
\|\theta_t\|_{L^2}^2+\|\theta\|_{L^\infty}^2(\|u\|_{L^\infty}^2\|\nabla \rho\|_{L^2}^2
+ \|\dv u\|_{L^2}^2)\Big)\\
 &\le  \delta \|\nabla u_t\|_{L^2}^2+ C_\delta\Big(
\|\theta_t\|_{L^2}^2 +
\|\theta\|_{H^1}^{2-\epsilon}\|\theta\|_{H^2}^{\epsilon}
( \|u\|_{H^1}^{2-\epsilon}\|u\|_{H^2}^{ \epsilon}+1)\Big),
\end{split}
\end{equation*}
for any $0<\epsilon<\frac{1}{2}$, which follows
\begin{equation}
\begin{split}
\int_0^t|K_1|ds  \le & \delta \|\nabla u_t\|_{L^2_t(L^2)}^2+
C_\delta \Big(1+
 \int_0^t\|\theta\|_{H^2}^\epsilon(\|u\|_{H^2}^\epsilon +1) ds\Big) \\
\le & \delta \|\nabla u_t\|_{L^2_t(L^2)}^2+ C_\delta ,
\end{split}
\end{equation}
by applying Lemmas \ref{lem2.5} and \ref{lemtheta}.

Next,
\begin{equation*}
\begin{split}
\int_0^t|K_2|ds \le & \int_0^t\int_\Omega \rho |u|(|u_t||\nabla u_t| +|\nabla u|^2 |u_t|  +
    |u|  |\nabla^2 u||u_t|+   |\nabla u||\nabla u_t|)dxds\\
\le & \int_0^t  \|u\|_{H^1}(\|u_t\|_{L^3}\|\nabla u_t\|_{L^2} +
\|\nabla u\|_{L^2}\|\nabla u\|_{H^1}\|u_t\|_{H^1}\\
&\qquad\quad +\|u\|_{H^1}\|u\|_{H^2}\|u_t\|_{H^1}+ \|\nabla u\|_{H^1}\|\nabla u_t\|_{L^2})ds\\
\le & \int_0^t (\|u_t\|_{L^2}^\frac{1}{2}\|u_t\|_{H^1}^\frac{3}{2}
+\|u\|_{H^2}\|u_t\|_{H^1} )ds\\
\le & \delta \|u_t\|_{L^2_t(H^1)}^2+
C_\delta,
\end{split}
\end{equation*}
and
\begin{equation*}
\begin{split}
\int_0^t|K_3|ds \le  \int_0^t\|\nabla
u\|_{L^2}\|u_t\|_{L^2}^\frac{1}{2}\|u_t\|_{H^1}^\frac{3}{2}ds \le
\delta \|u_t\|_{L^2_t(H^1)}^2+ C_\delta.
\end{split}
\end{equation*}

Note that by the boundary conditions \eqref{1.5} and the constraint
on viscosity coefficients $\mu,\lambda$, we have (see also
\cite{Lions98}, pp.76)
\begin{equation}\label{ut}
\begin{split}
\int_\Omega\Big(2\mu|D(u_t)|^2+ \lambda (\dv u_t)^2\Big)dx &\ge
2\mu\int_\Omega \Big(|D(u_t)|^2-\frac{1}{N} (\dv u_t)^2\Big)dx\\
&\ge \bar\mu \|u_t\|_{H^1}^2
\end{split}
\end{equation}
for some constant $\bar\mu>0$. Then we conclude from \eqref{2.17},
\eqref{ut} and the estimates for $K_1$ through $K_3$ that
\begin{equation}\label{22}
\begin{split}
\sup\limits_{0\le t\le T}\|u_t(t)\|_{L^2}^2 +
\|u_t\|_{L^2_t(H^1)}^2\le C .
\end{split}
\end{equation}

From  \eqref{1.2} and the elliptic regularity theory, we have
\begin{equation}\label{uh2}
\begin{split}
\|u\|_{L^\infty_t(H^2)} \le & C(\|u_t\|_{L^\infty_T(L^2)}
+ \|u\|_{L^\infty_t(H^1)}\|\nabla u\|_{L^\infty_t(L^3)}\\
&+ \|\nabla \theta\|_{L^\infty_T(L^2)}
+ \|\theta\|_{L^\infty_t(L^\infty)}\|\nabla \rho\|_{L^\infty_t(L^2)} )\\
\le & C(1+ \|u\|_{L^\infty_t(H^2)}^\frac{1}{2}+
\|\theta\|_{L^\infty_t(H^2)}^{ \epsilon}) ,
\end{split}
\end{equation}
for any $\epsilon>0$. Similarly, we derive from \eqref{1.3} that
\begin{equation}\label{th2}
\begin{split}
\|\theta\|_{L^\infty_t(H^2)} \le & C(\|\theta_t\|_{L^\infty_t(L^2)}+
\|u\|_{L^\infty_t(H^2)}).
\end{split}
\end{equation}

Substituting \eqref{uh2} into \eqref{th2} and choosing $\epsilon<1$,
we obtain
\begin{equation}\label{25}
\begin{split}
\|u\|_{L^\infty_t(H^2)} +\|\theta\|_{L^\infty_t(H^2)} \le &
C(1+\|\theta_t\|_{L^\infty_t(L^2)}).
\end{split}
\end{equation}

Next, we derive bounds for $\theta_t$ to close the desired energy estimates.
 Taking $\partial_t$ on both sides of the equation (\ref{1.3}),
then multiplying the resulting equation by $\theta_t$ in $L^2(\Omega )$,
we obtain
\begin{equation}\label{thetat}
\begin{split}
&\frac{1}{2}\frac{d}{dt}\int_\Omega\rho\theta_t^2dx +\kappa\int_\Omega|\nabla\theta_t|^2dx\\
& \quad = \int_\Omega R\rho \theta_t^2{\rm div}u dx+ \int_\Omega
R\rho_t\theta{\rm div}u\theta_tdx +\int_\Omega
R\rho \theta {\rm div}u_t\theta_tdx\\
&\qquad +\int_\Omega \big[ \mu (\nabla u+\nabla u^t):(\nabla u_t+\nabla u_t^t)
+2\lambda\dv u\,\dv u_t\big] \theta_tdx\\
&\qquad - \int_\Omega\rho_t u\cdot\nabla\theta\theta_tdx -\int_\Omega\rho
u_t\cdot\nabla\theta\theta_tdx-\int_\Omega\rho_t\theta_t^2dx\\
 & \quad  =:\sum\limits_{i=1}^7 L_i.
\end{split}
\end{equation}

We have to estimate each term on the right-hand side of
(\ref{thetat}). From \eqref{1.1} and Sobolev's
embedding theorem, we   get
\begin{equation*}
\begin{split}
\int_0^t |L_1|ds &\leq C \int _0^t\|\theta_t\|_{H^1}\| \theta_t\|_{L^2}\|{\rm div}u\|_{L^3}ds\\
&\leq \delta \|\theta_t\|_{L^2_t(H^1)}^2 +
C_\delta\|u\|_{L^\infty_t(H^1)}\|u\|_{L^\infty_t(H^2)} \|\theta_t\|_{L^2_t(L^2)}^2\\
&\le \delta \|\theta_t\|_{L^2_t(H^1)}^2 +
C_\delta\|u\|_{L^\infty_t(H^2)} , \label{2.18a}
\end{split}
\end{equation*}
\begin{equation*}
\begin{split}
\int_0^t |L_2|ds &\leq
\Big|\int_\Omega R(\rho {\rm div}\,u +\nabla \rho\cdot u)
\theta {\rm div}\,u \theta_t dx\Big|    \\
& \le   (\|\dv u\|_{L^\infty_t(L^2)}+ \|u\|_{L^\infty_{t}(H^1)}\|\nabla \rho\|_{L^\infty_t(L^2)})\\
&\qquad\times    \|\theta\|_{L^\infty_t(H^1)}\|\dv u\|_{L^2_t(H^1)}\|\theta_t\|_{L^2_t(H^1)}\\
& \le \delta\|\theta_t\|_{L^2_t(H^1)}^2 +C_\delta,  \label{2.18b}
\end{split}
\end{equation*}
 \begin{equation*}
\begin{split}
\int_0^t|L_3|ds & \leq
C\int_0^t\| \theta_t\|_{H^1}\|\theta\|_{H^1}\|{\rm div}u_t\|_{L^2}ds\qquad\\
&\le \delta \|\theta_t\|_{H^1}^2 + C_\delta, \label{2.18c}
\end{split}
\end{equation*}
\begin{equation*}
\begin{split}
\int_0^t|L_4|ds & \le  C\int_0^t\|\nabla u\|_{L^3}\|\nabla u_t\|_{L^2}\|\theta_t\|_{H^1}ds\hskip 4cm\\
&\le  \delta\|\theta_t\|_{L^2_t(H^1)}^2 + C_\delta \|u\|_{L^\infty_t(H^1)} \|u\|_{L^\infty_t(H^2)}\|\nabla u_t\|_{L^2_t(L^2)}^2\\
&\le  \delta\|\theta_t\|_{L^2_t(H^1)}^2 + C_\delta
 \|u\|_{L^\infty_t(H^2)} ,
\end{split}
\end{equation*}
\begin{equation*}
\begin{split}
\int_0^t|L_5|ds\le & C\int_0^t\int_\Omega (\rho |{\rm div}\,u|
+|u||\nabla\rho|)|u||\nabla\theta||\theta_t|dxds\\
\le & C\int_0^t(\|{\rm
div}u\|_{L^2}+\|\nabla\rho\|_{L^2}\|u\|_{H^1})
\|u\|_{H^1}\|\nabla\theta\|_{H^1}\|\theta_t\|_{H^1}ds\\
\le & \delta\|\theta_t\|_{L^2_t(H^1)}^2+ C_\delta\|\theta\|_{L^2_t(H^2)}^2 ds\\
\le & \delta\|\theta_t\|_{L^2_t(H^1)}^2+ C_\delta ,
\end{split}
\end{equation*}
\begin{equation*}
\begin{split}
\int_0^t|L_6|ds\le & C\int_0^t \|u_t\|_{H^1}\|\nabla\theta\|_{L^2}\|\theta_t\|_{L^3}ds\\
\le & \delta(\|u_t\|_{L^2_t(H^1)}^2+\|\theta_t\|_{L^2_t(H^1)}^2)
+ C_\delta\|\theta_t\|_{L^2_t(L^2)}^2,\\
\le & \delta\|\theta_t\|_{L^2_t(H^1)}^2+ C_\delta ,
\end{split}
\end{equation*}
\begin{equation*}
\begin{split}
\int_0^t|L_7|ds\le & C\int_0^t(\|\dv u\|_{L^2}+ \|u\|_{H^1}\|\nabla\rho\|_{L^2})\|\theta_t\|_{L^2}\|\theta_t\|_{H^1}ds\qquad\\
\le & \delta \|\theta_t\|_{L^2_t(H^1)}^2 + C_\delta.
\end{split}
\end{equation*}

Now, we integrate \eqref{thetat} and utilize the estimates for $L_1$ through $L_7$
with $\delta$ sufficiently small to conclude
\begin{equation}\label{tt}
\begin{split}
\|   \theta_t &(t)\|_{L^2}^2 + \|\theta_t\|_{L^2_t(H^1)}^2 \le C (1+
\|u\|_{L^\infty_t(H^2)} ),\quad 0\leq t\leq T.
\end{split}
\end{equation}
 As a consequence of \eqref{uh2}, \eqref{25} and \eqref{tt}, the current lemma is shown.
 \hfill$\square$
\vspace{2mm}

 Finally, in the next lemma we show the additional $L^q$
boundedness of the solution. The proof is exactly as in \cite{FJO},
however, we still reproduce it for the sake of completeness.
\begin{Lemma}\label{le2.5} Let $q$ be the same as in Theorem \ref{th1.1}. Then,
\begin{eqnarray}
\sup\limits_{0\leq t\leq T}\left(\|\rho_{t}(t)\|_{L^q}+\|\rho(t)\|_{W^{1,q}}
\right)\leq C,\label{2.22}\\
\int_{0}^{T}\Big(\|u(t)\|_{W^{2,q}}^{2}+\|\theta(t)\|_{W^{2,q}}^{2}\Big)dt\leq
C.\label{2.24}
\end{eqnarray}
\end{Lemma}
{\bf Proof.} Differentiating (\ref{1.1}) with respect to $x_{j}$ and multiplying
the resulting equation by $|\partial_{j}\rho|^{q-2}\partial_{j}\rho$ in
$L^2(\Omega )$, one deduces that
\begin{eqnarray*}
\frac{d}{dt}\int_\Omega |\nabla\rho|^q dx & \leq & C\int_\Omega
\big(
|\nabla u|\,|\nabla\rho|^q +|\rho|\,|\nabla\rho|^{q-1}|\nabla^{2}u|\big) dx  \\
& \leq & C\|\nabla u\|_{L^\infty}\|\nabla\rho\|^q_{L^q} +C\|\nabla^2
u\|_{L^q}\|\nabla\rho\|^{q-1}_{L^q},
\end{eqnarray*}
which gives
\begin{equation}
\begin{split}
 \sup_{0\leq t\leq T} \|\nabla\rho(t)\|_{L^q}
 &\leq  C\exp \Big(\int_0^t\|\nabla u(s)\|_{L^\infty}ds\Big)
 \Big(\|\nabla \rho_0\|_{L^q} + \int_0^t \|\nabla^2 u (s)\|_{L^q}ds\Big)\\
  &\le C(\sqrt{T})/\delta+ \delta \|\nabla^2 u \|_{L^2_t(L^q)} \label{rhoq}
 \end{split}
\end{equation}
 by Gronwall's inequality.

Using the regularity theory of elliptic equations again, we see that
\begin{eqnarray*}
\|u(t)\|_{W^{2,q}}&\leq&C\left(\|u_{t}\|_{L^q}+\|u\cdot\nabla
u\|_{L^q}+\|\nabla\rho\|_{L^q}+\|\nabla\theta\|_{L^q}\right)\\
&\leq&C\left(\|\nabla u_{t}\|_{L^{2}}+\|u\|_{L^{\infty}}\|\nabla
u\|_{L^q}+\|\nabla\rho\|_{L^q}+\|\theta\|_{H^{2}}\right)\\
&\leq&C\left(\|\nabla
u_{t}\|_{L^{2}}+\|u\|_{H^{2}}^{2}+\|\nabla\rho\|_{L^q}+\|\theta\|_{H^{2}}\right).
\end{eqnarray*}
If we integrate the above inequality over $(0,T)$ and make use of
(\ref{rhoq}) as well as the estimates we have proved so far, we obtain
\begin{equation}
\int_{0}^{T}\|u(t)\|_{W^{2,q}}^{2}dt\leq C, \label{last}
\end{equation}
and thus, from \eqref{rhoq} one gets
$$\sup\limits_{0\le t\le T}\|\rho(t)\|_{W^{1,q}}\le
C.$$

Since $\rho_{t}=-u\nabla\rho-\rho\dv u$, we also have
$$\|\rho_{t}(t)\|_{L^q}\leq\|u\|_{L^{\infty}}\|\nabla\rho\|_{L^q}
+\|\rho\|_{L^{\infty}}\|\dv u\|_{L^q}\leq C.$$  Then the boundedness
of $\theta$ in $L^2(0,T;W^{2,q})$ follows from (\ref{1.3}),
(\ref{last}) and the above inequality. The proof is finished.
\hfill$\square$
\subsection{Conclusions.} By virtue of all the above energy
estimates, we obtain the bounds of the norms of $(\rho,u,\theta)$ in
$[0,T]\times\Omega$ in the sense of \eqref{1.7} for any $T<T^*$.
These bounds depend only on $\Omega$, the initial data, and
continuously on $T^*$ (in fact, the bounds depend on $T^*$ either
polynomially or exponentially!). Thus, we can take
$(\rho,u,\theta,\rho_t,  u_t, \theta_t)|_{t=T}$, redefined if
necessary, as the initial data at $t=T$ and apply Proposition
\ref{pr1.1} to extend the solution to $t=T+T_1$.

If $T+T_1>T^*$, then it contradicts the maximality of $T^*$.
Otherwise, we can continue to extend the solution by
taking the values of the solution at $t=T+T_1$ as initial data
again. Since the a priori estimates are independent of any $t<T^*$,
the solution can be extended to $t=T+2T_1$. Here, we remark that
by applying Proposition \ref{pr1.1}, the solution can be extended
from $t=T+T_1$ to $t=T+2T_1$, since the local existence interval
depends only on the initial data which, in our case, are bounded
in any time interval $[0,\overline{T}]$ with a bound depending on
$\overline{T}$ only. Utilizing Proposition
\ref{pr1.1} repeatedly, there must exist a positive integer $m$, such
that $T+mT_1>T^*$. This also leads to the contradiction to the
maximality of $T^*$. Therefore, the assumption (\ref{2.1}) does not hold.
This completes the proof of Theorem \ref{th1.1}.
\\[5mm]
{\bf Acknowledgements.} The research of Ou is partially supported by the China Postdoctoral
Science Foundation (Grant No. 20090450333).
Jiang is supported by the National Basic Research Program (Grant No. 2005CB321700)
and NSFC (Grant No. 40890154).

\end{document}